\documentclass{amsart}

\usepackage{amssymb}

\allowdisplaybreaks

\newtheorem{theorem}{Theorem}[section]
\newtheorem{lemma}[theorem]{Lemma}
\newtheorem{theoA}{Theorem}
\newtheorem{theoB}{Theorem}
\newtheorem{corC}{Corollary}

\newtheorem{remark}[theorem]{Remark}

\newtheorem{problem}[theorem]{Problem}

\theoremstyle{definition}

\newcommand{\N}{\mathbb{N}}

\newcommand{\R}{\mathbb{R}}
\newcommand{\C}{\mathbb{C}}

\newcommand{\summ}{\sum\nolimits}

\newcommand{\dem}{\noindent {\bf Proof. }}
\newcommand{\fin}{\hspace*{\fill} $\square$ \vskip0.2cm}

\begin{document}

\title[Weak Burkholder inequality]
{Weak type estimates associated to \\ Burkholder's martingale
inequality}

\author[Javier Parcet]
{Javier Parcet}

\address{Universidad Aut\'{o}noma de Madrid}

\email{javier.parcet@uam.es}

\footnote{Partially supported by the Project MTM2004/00678,
Spain.} 

\begin{abstract}
Given a probability space $(\Omega, \mathsf{A}, \mu)$, let
$\mathsf{A}_1, \mathsf{A}_2, \ldots$ be a filtration of
$\sigma$-subalgebras of $\mathsf{A}$ and let $\mathsf{E}_1,
\mathsf{E}_2, \ldots$ denote the corresponding family of
conditional expectations. Given a martingale $f = (f_1, f_2,
\ldots)$ adapted to this filtration and bounded in $L_p(\Omega)$
for some $2 \le p < \infty$, Burkholder's inequality claims that
$$\|f\|_{L_p(\Omega)} \sim_{\mathrm{c}_p} \Big\| \Big(
\sum_{k=1}^\infty \mathsf{E}_{k-1}(|df_k|^2) \Big)^{1/2}
\Big\|_{L_{p}(\Omega)} + \Big( \sum_{k=1}^\infty \|df_k\|_p^p
\Big)^{1/p}.$$ Motivated by quantum probability, Junge and Xu
recently extended this result to the range $1 < p < 2$. In this
paper we study Burkholder's inequality for $p=1$, for which the
techniques (as we shall explain) must be different. Quite
surprisingly, we obtain two non-equivalent estimates which play
the role of the weak type $(1,1)$ analog of Burkholder's
inequality. As application, we obtain new properties of Davis
decomposition for martingales.
\end{abstract}

\maketitle

\section*{Introduction and Main Results}

\renewcommand{\theequation}{$\mathsf{B}_{p}^+$}
\addtocounter{equation}{-1}

Sums of independent random variables and martingale inequalities
are nowadays powerful tools in classical harmonic analysis. Mostly
in the 70's and 80's, the works of Bourgain, Burkholder, Davis,
Gundy, Pisier, Rosenthal and many others illustrated a fruitful
interaction between these subjects and Calder{\'o}n-Zygmund theory of
singular integrals as well as Littlewood-Paley theory. Our first
motivation in this paper comes from a fundamental result due to
Burkholder \cite{B,BG} which can be stated as follows. Given a
probability space $(\Omega, \mathsf{A}, \mu)$, let $\mathsf{A}_1,
\mathsf{A}_2, \ldots$ be a filtration of $\sigma$-subalgebras of
$\mathsf{A}$ and let $\mathsf{E}_1, \mathsf{E}_2, \ldots$ denote
the corresponding family of conditional expectations. Given $2 \le
p < \infty$ and a martingale $f = (f_1, f_2, \ldots)$ adapted to
this filtration and bounded in $L_p(\Omega)$, we have
\begin{equation} \label{Eq-Burk-p>2}
\|f\|_{p} \sim_{\mathrm{c}_p} \Big\| \Big( \sum_{k=1}^\infty
\mathsf{E}_{k-1}(|df_k|^2) \Big)^{1/2} \Big\|_{p} + \Big\| \Big(
\sum_{k=1}^\infty |df_k|^p \Big)^{1/p} \Big\|_{L_{p}(\Omega)}.
\end{equation}
The first term on the right is called the \emph{conditional square
function} of $f$, while the second is the \emph{$p$-variation} of
$f$. The optimal growth of the equivalence constant is given by
$\mathrm{c}_p \sim p / \log p$ as $p \to \infty$ and
\eqref{Eq-Burk-p>2} fails on $L_\infty(\Omega)$, we refer to the
papers \cite{H,JSZ} for more details. Apart from the relation with
harmonic analysis, Burkholder's inequality has deep implications
in the geometry of Banach spaces. Let us mention for instance
Maurey/Pisier's theory of type and cotype or the isomorphism and
embedding theory of $L_p$ spaces via $p$-stable processes.

\vskip5pt

These assertions are justified by the following observations.
First, Rosenthal's inequality \cite{Ro} appears as the particular
case where the sequence $df_1, df_2, \ldots$ is given by a family
of independent mean-zero random variables $df_k = \xi_k$. In this
case we have $\mathsf{E}_{k-1}(|df_k|^2) = \|\xi_k\|_2^2$ and
deduce Rosenthal's inequality
\renewcommand{\theequation}{\arabic{equation}}
\renewcommand{\theequation}{$\mathsf{R}_p$}
\addtocounter{equation}{-1}
\begin{equation} \label{Eq-Rosenthal}
\Big\| \sum_{k=1}^\infty \xi_k \Big\|_p \sim_{\mathrm{c}_p} \Big(
\sum_{k=1}^\infty \|\xi_k\|_2^2 \Big)^{1/2} + \Big(
\sum_{k=1}^\infty \|\xi_k\|_p^p \Big)^{1/p}.
\end{equation}
Moreover, we can go further and take $\xi_k = \lambda_k
\varepsilon_k$ with $\lambda_k \in \C$ and $\varepsilon_1,
\varepsilon_2, \varepsilon_3, \ldots$ independent Bernoulli random
variables equidistributed in $\pm 1$, the reader can think for
instance in the sequence of Rademacher functions on the unit
interval. In this case, the two terms on the right of
\eqref{Eq-Rosenthal} collapse into the first one and we recover
the classical Khintchine inequalities for $2 \le p < \infty$
\renewcommand{\theequation}{\arabic{equation}}
\renewcommand{\theequation}{$\mathsf{K}_p$}
\addtocounter{equation}{-1}
\begin{equation} \label{Eq-Khintchine}
\Big\| \sum_{k=1}^\infty \lambda_k \varepsilon_k \Big\|_p
\sim_{\mathrm{c}_p} \Big( \sum_{k=1}^\infty |\lambda_k|^2
\Big)^{1/2}.
\end{equation}

Our second motivation comes from the noncommutative analogues of
the results mentioned so far. Roughly speaking, we replace
functions by operators (this process is known as
\emph{quantization}) and study noncommutative generalizations of
the classical results. In this setting, the main objects are
noncommutative $L_p$ spaces constructed over von Neumann algebras
\cite{PX2} and noncommutative martingales \cite{X}. The theory of
noncommutative martingale inequalities (a subfield of quantum
probability) has reached an espectacular development after
Pisier/Xu's seminal paper \cite{PX1}. Indeed, it can be said that
almost every classical martingale inequality has been successfully
transferred to the noncommutative setting. We find noncommutative
analogues of the Burkholder-Gundy inequalities \cite{PX1}, Doob's
maximal inequality \cite{J1} and weak type $(1,1)$ estimates for
martingale transforms \cite{R}. Burkholder's inequality was
finally obtained by Junge and Xu in \cite{JX}.

\vskip5pt

The new insight provided by the noncommutative formulation led
Junge and Xu to extend in \cite{JX} Burkholder's inequality to the
range $1 < p \le 2$. This extension was new even in the
commutative case and can be explained (we only consider here
commutative random variables) as follows. The right hand side of
\eqref{Eq-Burk-p>2} can be understood as the norm in the
intersection of two Banach spaces, respectively called
\emph{conditional} and \emph{diagonal} Hardy spaces of
martingales. Thus, it is natural to guess that in the dual
formulation of \eqref{Eq-Burk-p>2}, we will find a \emph{sum} of
the dual Hardy spaces. Then, recalling the definition of the norm
on a sum of Banach spaces, Burkholder's inequality for $1 < p \le
2$ reads as follows
\renewcommand{\theequation}{\arabic{equation}}
\renewcommand{\theequation}{$\mathsf{B}_{p}^-$}
\addtocounter{equation}{-1}
\begin{equation} \label{Eq-Burk-p<2}
\|f\|_{p} \, \sim_{\mathrm{c}_p} \, \inf_{f=g+h} \left\{ \Big\|
\Big( \sum_{k=1}^\infty \mathsf{E}_{k-1}(|dg_k|^2) \Big)^{1/2}
\Big\|_{p} + \Big\| \Big(\sum_{k=1}^\infty |dh_k|^p \Big)^{1/p}
\Big\|_p \right\},
\end{equation}
where the infimum runs over all possible decompositions of $f$ as
a sum $f=g+h$ of two martingales. Note that the right expressions
for ($\mathsf{B}_2^+$) and ($\mathsf{B}_2^-$) are equivalent so
that the quadratic case \emph{explains} the transition from
intersections to sums. Let us mention that this is a typical
phenomenon in the noncommutative setting which also appears for
instance in the noncommutative Khintchine inequalities \cite{LuP}
or the noncommutative Burkholder-Gundy inequalities \cite{PX1}. In
particular, Junge/Xu's paper \cite{JX} illustrated how a
noncommutative problem can give some light in its commutative
counterpart! We think this is also the case in this paper.

\vskip5pt

The problem of determining the behavior of Burkholder's inequality
on $L_1(\Omega)$ naturally came out after Junge/Xu's extension. As
it was pointed out in \cite{JX}, the upper estimate holds with an
absolute constant $\mathrm{c}$
$$\|f\|_1 \, \le \mathrm{c} \inf_{f=g+h} \left\{ \Big\| \Big(
\sum_{k=1}^\infty \mathsf{E}_{k-1}(|dg_k|^2) \Big)^{1/2} \Big\|_1
+ \Big\| \sum_{k=1}^\infty |dh_k| \Big\|_1 \right\}.$$ Keeping the
notation, our first weak type estimate reads as follows.
\renewcommand{\theequation}{\arabic{equation}}
\renewcommand{\theequation}{$\mathsf{WB}_1$}
\addtocounter{equation}{-1}

\begin{theoA}
Let $f = (f_1, f_2, \ldots)$ be a bounded martingale in
$L_1(\Omega)$. Then, we can decompose $f$ as a sum $f= g+h$ of two
martingales adapted to the same filtration and satisfying the
following inequality with an absolute constant $\mathrm{c}$
\begin{equation} \label{Eq-Weak-Burkholder}
\Big\| \Big( \sum_{k=1}^\infty \mathsf{E}_{k-1}(|dg_k|^2)
\Big)^{1/2} \Big\|_{1,\infty} + \Big\| \sum_{k=1}^\infty |dh_k|
\Big\|_{1,\infty} \le \mathrm{c} \, \|f\|_1.
\end{equation}
\end{theoA}

The first relevant difficulty in proving Theorem A resides on the
fact that we need to guess the \emph{right} decomposition of $f$
before estimating the associated norms on the corresponding Hardy
spaces. Note that this problem was avoided in \cite{JX} by using a
duality argument, which is not available in our case. The solution
to this problem turns out to be very elegant. Indeed, quite
surprisingly the right decomposition is given by the classical
Davis decomposition \cite{D}. Moreover, the proof becomes quite
involved since we combine Davis and Gundy decompositions to
estimate the diagonal term. It is also worthy of mention that
Theorem A provides an \emph{improvement} of Davis decomposition.
The reader is referred to the last section of this paper for the
details.

\vskip5pt

The first term on the left of \eqref{Eq-Weak-Burkholder} is
clearly the natural \emph{weak} analog of the first term on the
right of \eqref{Eq-Burk-p<2}. However, the second term on the left
of \eqref{Eq-Weak-Burkholder} is only one possible interpretation
for the corresponding $L_p$ term. Indeed, we have chosen the
$L_{1,\infty}$ norm of the $1$-variation of $h$. In other words,
the norm of the martingale difference sequence $dh$ in
$L_{1,\infty}(\Omega; \ell_1)$. This choice is motivated by the
fact that the $L_p$ term is exactly the norm of $dh$ in
$L_p(\Omega; \ell_p)$. Another interpretation arises after
rewriting $L_p(\Omega; \ell_p)$ as the scalar-valued space
$L_p(\Omega_{\oplus \infty})$ where the associated measure is
given by $$\mu_{\oplus \infty} \Big( \bigoplus_{k \ge 1}
\mathrm{A}_k \Big) = \sum_{k \ge 1} \mu(\mathrm{A}_k).$$ In this
case, the weak analog of the $p$-variation is given by
$$\Big\| \sum_{k=1}^{\infty} \delta_k \otimes dh_k
\Big\|_{L_{1,\infty}(\Omega_{\oplus \infty})} = \ \sup_{\lambda
> 0} \ \lambda \sum_{k=1}^{\infty} \mu \Big\{ |dh_k|
> \lambda \Big\},$$ where $(\delta_k)_{k \ge 1}$ denotes the
canonical unit vector basis. At this point, it is worthy of
mention that the norms of $L_{1,\infty}(\Omega;\ell_1)$ and
$L_{1,\infty}(\Omega_{\oplus \infty})$ are not equivalent nor even
comparable. Indeed, taking $\varphi_k = \chi_{[0,1/k]}$ and $\xi_k
= \frac1k \chi_{[0,1]}$, it is easy to check that $$\sup_{\lambda
> 0} \, \lambda \mu \Big\{ \sum_{k=1}^m \varphi_k > \lambda \Big\} \sim
1 << \log m \sim \sup_{\lambda
> 0} \, \lambda \sum_{k=1}^m \mu \Big\{ \varphi_k > \lambda \Big\},$$
$$\sup_{\lambda > 0} \, \lambda \mu \Big\{ \sum_{k=1}^m \xi_k >
\lambda \Big\} \sim \log m >> 1 \sim \sup_{\lambda
> 0} \, \lambda \sum_{k=1}^m \mu \Big\{ \xi_k > \lambda \Big\}.$$
Therefore, it makes sense to study the weak type estimate
associated to this new interpretation of the $p$-variation. In
order to state our second weak type estimate, we need to recall
the notion of a regular filtration. The filtration $\mathsf{A}_1,
\mathsf{A}_2, \mathsf{A}_3, \ldots$ of $\sigma$-subalgebras of
$\mathsf{A}$ is called \emph{$\mathrm{k}$-regular} for some
constant $\mathrm{k} > 1$ if every non-negative martingale $f=
(f_1, f_2, \ldots)$ adapted to this filtration satisfies $$f_n \le
\mathrm{k} f_{n-1}.$$ Examples of regular martingales arise from
the filtrations generated by bounded Vilenkin systems. In
particular, the dyadic martingales are the most well-known
examples of this kind, see \cite{W} for more on this topic.
\renewcommand{\theequation}{\arabic{equation}}
\renewcommand{\theequation}{$\mathsf{WB}_2$}
\addtocounter{equation}{-1}

\begin{theoB}
Let $f = (f_1, f_2, \ldots)$ be a bounded martingale in
$L_1(\Omega)$ adapted to a $\mathrm{k}$-regular filtration. Then,
we can decompose $f$ as a sum $f= g+h$ of two martingales adapted
to the same filtration and satisfying the following inequality
with an absolute constant $\mathrm{c}$
\begin{equation} \label{Eq-Weak-Burkholder2}
\Big\| \Big( \sum_{k=1}^\infty \mathsf{E}_{k-1}(|dg_k|^2)
\Big)^{1/2} \Big\|_{L_{1,\infty}(\Omega)} + \Big\|
\sum_{k=1}^\infty \delta_k \otimes dh_k
\Big\|_{L_{1,\infty}(\Omega_{\oplus \infty})} \le \mathrm{c}
\mathrm{k} \, \|f\|_1.
\end{equation}
\end{theoB}

The notion of $\mathrm{k}$-regular filtration (equivalently that
of previsible martingale) is necessary to formulate many
martingale inequalities, see e.g. \cite{B,BG} or Chapter 2 in
\cite{W}. However, it is still unclear whether or not the
$\mathrm{k}$-regularity assumption in Theorem B is necessary. On
one side, in view of some similar results in \cite{BG}, it seems a
natural condition. However, the proof we present here (see
Paragraph \ref{Paragraph3.4} for a much simpler but less
interesting one) gives some evidences that Theorem B might hold
for general martingales with an absolute constant.

\vskip5pt

Nevertheless, even in the present form, Theorem B presents some
advantages with respect to Theorem A. First, as we shall explain
in the last section of this paper, it is much simpler to reprove
Burkholder's inequality (via real interpolation and duality)
starting from Theorem B. Second, our \emph{more elaborated} proof
of Theorem B goes further and gives rise to the result below.

\begin{corC}
Let $f = (f_1, f_2, \ldots)$ be a bounded martingale in
$L_1(\Omega)$. Then, we can decompose each $f_n$ as a sum $f_n =
g_n + h_n$ of two functions $($non-necessarily martingales$)$
adapted to the same filtration and satisfying the following
inequality with an absolute constant $\mathrm{c}$ $$\Big\| \Big(
\sum_{k=1}^\infty \mathsf{E}_{k-1}(|dg_k|^2) \Big)^{1/2}
\Big\|_{L_{1,\infty}(\Omega)} + \Big\| \sum_{k=1}^\infty \delta_k
\otimes dh_k \Big\|_{L_{1,\infty}(\Omega_{\oplus \infty})} \le
\mathrm{c} \, \|f\|_1.$$
\end{corC}

At the time of this writing, Randrianantoanina independently
obtained in \cite{R4} the noncommutative analogue of the result
above. However, Theorems A and B have not been considered there.
The noncommutative form of Corollary C has been applied to obtain
the optimal constants in Junge/Xu's noncommutative Burkholder
inequality. The reader familiar with the noncommutative setting
will recognize the similarities between both papers. In any case,
we have decided to add a paragraph at the end explaining how both
arguments are related. More concretely, although this is not
mentioned in Randrianantoanina's paper, we shall explain how Davis
decomposition appears (in a very indirect form) in \cite{R4}.

\vskip5pt

\noindent \textbf{Acknowledgements.} I have received interesting
suggestions and comments from Teresa Mart{\'\i}nez, Marius Junge and
Fernando Soria. I am specially imdebted to Michael Cwikel for
disproving the isomorphism \eqref{Eq-RealInt} below.

\section{Martingale Decompositions}
\label{Section1}

\renewcommand{\theequation}{\arabic{equation}}

Leu us fix once and for all a probability space $(\Omega,
\mathsf{A}, \mu)$ and a filtration $\mathsf{A}_1, \mathsf{A}_2,
\ldots$ of $\sigma$-subalgebras of $\mathsf{A}$ with corresponding
conditional expectations $\mathsf{E}_1, \mathsf{E}_2, \ldots$
Davis decomposition is a fundamental tool in the theory of
martingale inequalities and it appeared for the first time in
\cite{D}, where Davis applied it to prove his well-known theorem
on the equivalence in $L_1(\Omega)$ between the martingale square
function and Doob's maximal function
$$\|f^*\|_1 \sim_{\mathrm{c}} \Big\| \Big( \sum_{k=1}^\infty
|df_k|^2 \Big)^{1/2} \Big\|_1.$$ Considering the truncated maximal
functions $$f_n^*(w) = \sup_{1 \le k \le n} |f_k(w)|,$$ we
formulate Davis decomposition $f=g+h$ by defining the differences
\begin{eqnarray*}
dg_k & = & df_k \chi_{\{f_k^* < 2 f_{k-1}^*\}} - \mathsf{E}_{k-1}
\big( df_k \chi_{\{f_k^* < 2 f_{k-1}^*\}} \big), \\ dh_k & = &
df_k \chi_{\{f_k^* \ge 2 f_{k-1}^*\}} - \mathsf{E}_{k-1} \big(
df_k \chi_{\{f_k^* \ge 2 f_{k-1}^*\}} \big).
\end{eqnarray*}
It is clear that $dg_k$ and $dh_k$ are martingale differences so
that $g$ and $h$ become martingales adapted to the filtration
$\mathsf{A}_1, \mathsf{A}_2, \ldots$ The properties stated in
\cite{D} and which appear in the literature are the following
\begin{equation} \label{Eq-Properties-Davis}
|dg_k| \le 8 f_{k-1}^* \quad \mbox{and} \quad \Big\|
\sum_{k=1}^\infty |dh_k| \Big\|_p \le (4+4p) \, \|f^*\|_p
\end{equation}
for $1 \le p < \infty$, see e.g. \cite{W} for an estimate of $g$
in the norm of the space of predictable martingales. The proof of
these properties is rather simple, in contrast with their weak
type analogs which arise from Theorem A, see Section
\ref{Section4}.

\vskip5pt

Let us now describe Gundy's decomposition. Let $f=(f_1, f_2,
\ldots)$ be a martingale on $(\Omega, \mathsf{A},\mu)$ relative to
the filtration fixed above that is bounded in $L_1(\Omega)$. Let
$\lambda$ be a positive real number. Then we define the
martingales $\alpha$, $\beta$ and $\gamma$ by their martingale
differences
\begin{eqnarray}
\nonumber d \alpha_k & = & df_k \chi_{\{f_{k-1}^* > \lambda\}}, \\
\label{Eq-Gundy-Dec} d \beta_k & = & df_k \chi_{\{f_k^* \le
\lambda\}} -
\mathsf{E}_{k-1} \big( df_k \chi_{\{f_k^* \le \lambda\}} \big), \\
\nonumber d \gamma_k & = & df_k \chi_{\{f_{k-1}^* \le \lambda <
f_k^*\}} - \mathsf{E}_{k-1} \big( df_k \chi_{\{f_{k-1}^* \le
\lambda < f_k^*\}} \big).
\end{eqnarray}
Again, these are clearly martingale differences with sum $df_k$
and thus we get a decomposition $f = \alpha + \beta + \gamma$ into
three martingales. This decomposition is in fact due to Burkholder
\cite{B} and is simpler than the one originally formulated by
Gundy \cite{G}. Indeed, the decomposition stated above uses only
one stopping time while the one originally formulated by Gundy
needs two stopping times. I learned this simpler (but weaker, see
below) decomposition from Narcisse Randrianantoanina. The
following are the properties satisfied by the given decomposition
\begin{itemize}
\item[i)] The martingale $\alpha$ satisfies $$\lambda \, \mu
\Big\{ \sum_{k=1}^\infty |d \alpha_k| > 0 \Big\} \le \mathrm{c}
\|f\|_1.$$ \item[ii)] The martingale $\beta$ satisfies \vskip1pt
$$\|\beta\|_1 \le \mathrm{c} \|f\|_1 \quad \mbox{and} \quad
\frac{1}{\lambda} \, \|\beta\|_2^2 \le \mathrm{c} \|f\|_1.$$
\vskip1pt \null \item[iii)] The martingale $\gamma$ satisfies
$$\sum_{k=1}^{\infty} \|d \gamma_k\|_1 \le \mathrm{c} \|f\|_1.$$
\end{itemize}

Gundy's original decomposition, paying the price of using two
stopping times, obtains the additional estimate $\|\beta\|_\infty
\le \mathrm{c} \lambda$. In particular, we can control any $L_p$
norm of $\beta$ by means of H{\"o}lder's inequality. Nevertheless, we
shall only need the $L_2$ estimate in this paper. Gundy's
decomposition theorem plays a central role in classical martingale
theory and it can be regarded as a probabilistic counterpart of
the well-known Calder\'{o}n-Zygmund decomposition for integrable
functions in harmonic analysis, see \cite{GA,G} for further
details.

\section{Proof of Theorem A}
\label{Section2}

We first observe that we may assume without lost of generality
that $f$ is a positive martingale. Indeed, otherwise we can always
decompose $f_n$ into a linear combination of four positive
functions $$f_n = \big( f_n^{(1)} - f_n^{(2)} \big) + i \big(
f_n^{(3)} - f_n^{(4)} \big).$$ According to a classical result due
to Krickeberg, it turns out that this provides a martingale
decomposition of $f$ into four positive martingales. Therefore,
since we have the inequality
$$\sum_{k=1}^4 \|f_n^{(k)}\|_1 \le 2 \, \|f_n\|_1$$ and the
expressions on the left of \eqref{Eq-Weak-Burkholder} clearly
satisfy a quasi-triangle inequality, we may assume that $f$ is
positive. On the other hand, as we have anticipated in the
Introduction, the right decomposition to prove the weak Burkholder
inequality is Davis decomposition
\begin{eqnarray*}
dg_k & = & df_k \chi_{\{f_k^* < 2 f_{k-1}^*\}} - \mathsf{E}_{k-1}
\big( df_k \chi_{\{f_k^* < 2 f_{k-1}^*\}} \big), \\ dh_k & = &
df_k \chi_{\{f_k^* \ge 2 f_{k-1}^*\}} - \mathsf{E}_{k-1} \big(
df_k \chi_{\{f_k^* \ge 2 f_{k-1}^*\}} \big).
\end{eqnarray*}
In what follows, $\mathrm{c}$ might have different values from one
instance to another.

\subsection{Step 1:} \textbf{Proof of the estimate} $$\Big\| \Big(
\sum_{k=1}^{\infty} \mathsf{E}_{k-1} (|dg_k|^2) \Big)^{1/2}
\Big\|_{1,\infty} \le \mathrm{c} \, \|f\|_1.$$ Taking $\tau_k =
df_k \chi_{\{ f_k^* < 2 f_{k-1}^* \}}$, we have
$$\mathsf{E}_{k-1}(|dg_k|^2) = \mathsf{E}_{k-1} \Big( |\tau_k|^2 +
\big| \mathsf{E}_{k-1}(\tau_k) \big|^2 - 2 \, \mathrm{Re}
[\overline{\tau_k} \, \mathsf{E}_{k-1}(\tau_k)] \Big) \le
\mathsf{E}_{k-1} (|\tau_k|^2).$$ Thus, we may replace $dg_k$ by
$\tau_k$ and defining the function $$\Phi = \sum_{k=1}^\infty
\mathsf{E}_{k-1} \big( |df_k|^2 \chi_{\{f_k^* < 2 f_{k-1}^*\}}
\big),$$ it suffices to prove that $\lambda \mu \big( \Phi >
\lambda^2 \big) \le \mathrm{c} \|f\|_1$ for all $\lambda > 0$. For
fixed $\lambda$ we have $$\lambda \mu \big( \Phi > \lambda^2 \big)
\le \lambda \mu \big( f^* > \lambda \big) + \lambda \mu \big(
\chi_{\{f^* \le \lambda\}} \Phi > \lambda^2 / 2 \big).$$
Therefore, since $\mathsf{E}_0 = \mathsf{E}_1$ and $\chi_{\{f^*
\le \lambda\}} \le \chi_{\{f_{k-1}^* \le \lambda\}}$ for $k \ge
2$, we have $$\lambda \mu \big( \Phi > \lambda^2 \big) \le \lambda
\mu \big( f^* > \lambda \big) + \lambda \mu \Big(
\sum_{k=1}^\infty \mathsf{E}_{k-1} \big( |df_k|^2 \chi_{\{f_k^*
\le 2 \lambda\}} \big) > \lambda^2 / 2 \Big).$$ According to
Doob's maximal inequality, the first term is controlled by
$\mathrm{c} \|f\|_1$. On the other hand, in order to estimate the
second term, we need to decompose it into two pieces. We use
Chebychev's inequality
\begin{eqnarray*}
\lefteqn{\lambda \mu \Big( \sum_{k=1}^\infty \mathsf{E}_{k-1}
\big( |df_k|^2 \chi_{\{f_k^* \le 2 \lambda\}} \big)
> \lambda^2 / 2 \Big)} \\ & \le & \frac{2}{\lambda} \Big\|
\sum_{k=1}^\infty \mathsf{E}_{k-1} \big( |df_k|^2 \chi_{\{f_k^*
\le 2 \lambda\}} \big) \Big\|_1 \\ & = & \frac{2}{\lambda}
\sum_{k=1}^\infty \big\| df_k \chi_{\{f_k^* \le 2 \lambda\}}
\big\|_2^2 \\ & \le & \frac{4}{\lambda} \sum_{k=1}^\infty \big\|
f_k \chi_{\{f_k^* \le 2 \lambda\}} - f_{k-1} \chi_{\{f_{k-1}^* \le
2 \lambda\}} \big\|_2^2 \\ & + & \frac{4}{\lambda}
\sum_{k=2}^\infty \big\| f_{k-1} \chi_{\{f_{k-1}^* \le 2
\lambda\}} - f_{k-1} \chi_{\{f_{k}^* \le 2 \lambda\}} \big\|_2^2.
\end{eqnarray*}
Denoting by $\mathrm{A}$ and $\mathrm{B}$ the two terms on the
right, we have
\begin{eqnarray*}
\mathrm{A} & = & \frac{4}{\lambda} \sum_{k=1}^\infty \int_\Omega
f_k^2 \chi_{\{f_k^* \le 2 \lambda\}} - f_{k-1}^2 \chi_{\{f_{k-1}^*
\le 2 \lambda\}} \, d\mu \\ & + & \frac{8}{\lambda}
\sum_{k=2}^\infty \int_\Omega f_{k-1} \chi_{\{f_{k-1}^* \le 2
\lambda\}} \big( f_{k-1} - f_k \chi_{\{f_k^* \le 2 \lambda\}}
\big) \, d \mu \\ & = & \frac{4}{\lambda} \, \lim_{n \to \infty}
\int_\Omega f_n^2 \chi_{\{f_n^* \le 2 \lambda\}} \, d \mu
\\ & + & \frac{8}{\lambda}
\sum_{k=2}^\infty \int_\Omega f_{k-1} \chi_{\{f_{k-1}^* \le 2
\lambda\}} \big( f_{k-1} \chi_{\{f_{k-1}^* \le 2 \lambda\}} -
\mathsf{E}_{k-1} (f_k \chi_{\{f_k^* \le 2 \lambda\}}) \big) \, d
\mu \\ & \le & \frac{4}{\lambda} \, \lim_{n \to \infty} \|f_n\|_1
\, \big\| f_n \chi_{\{f_n^* \le 2 \lambda\}} \big\|_\infty
\\ & + & \frac{8}{\lambda}
\sum_{k=2}^\infty \big\| f_{k-1} \chi_{\{f_{k-1}^* \le 2
\lambda\}} \big\|_\infty \int_\Omega \big| f_{k-1}
\chi_{\{f_{k-1}^* \le 2 \lambda\}} - \mathsf{E}_{k-1} (f_k
\chi_{\{f_k^* \le 2 \lambda\}}) \big| \, d \mu.
\end{eqnarray*}
This gives $$\mathrm{A} \le 8 \, \|f\|_1 + 16 \sum_{k=2}^\infty
\int_\Omega \big| f_{k-1} \chi_{\{f_{k-1}^* \le 2 \lambda\}} -
\mathsf{E}_{k-1} (f_k \chi_{\{f_k^* \le 2 \lambda\}}) \big| \, d
\mu.$$ However, we may drop the modulus in the integral since
$$\mathsf{E}_{k-1} (f_k \chi_{\{f_k^* \le 2 \lambda\}}) \le
\mathsf{E}_{k-1} (f_k \chi_{\{f_{k-1}^* \le 2 \lambda\}}) =
f_{k-1} \chi_{\{f_{k-1}^* \le 2 \lambda\}}.$$ In that case,
$\mathsf{E}_{k-1}$ disappears and we obtain a telescopic sum of
decreasing terms $$\mathrm{A} \le 8 \, \|f\|_1 + 16
\sum_{k=2}^\infty \int_\Omega f_{k-1} \chi_{\{f_{k-1}^* \le 2
\lambda\}} \, d \mu - \int_\Omega f_k \chi_{\{f_k^* \le 2
\lambda\}} \, d \mu \le 24 \, \|f\|_1.$$ It remains to estimate
the term B. This term is much easier to handle
\begin{eqnarray*}
\mathrm{B} & = & \frac{4}{\lambda} \sum_{k=2}^\infty \big\|
f_{k-1} \chi_{\{f_{k-1}^* \le 2 \lambda < f_{k}^*\}} \big\|_2^2
\\ & \le & \frac{4}{\lambda} \sum_{k=2}^\infty \big\| \chi_{\{f_{k-1}^*
\le 2 \lambda < f_{k}^*\}} \big\|_1 \, \big\| f_{k-1}^2
\chi_{\{f_{k-1}^* \le 2 \lambda < f_{k}^*\}} \big\|_\infty \\ &
\le & 16 \lambda \sum_{k=2}^\infty \int_\Omega \chi_{\{f_{k-1}^*
\le 2\lambda\}} - \chi_{\{f_{k}^* \le 2\lambda\}} \, d \mu \le 16
\lambda \, \mu \big( f^* > 2 \lambda \big) \le 8 \, \|f\|_1.
\end{eqnarray*}

\subsection{Step 2:} \textbf{Proof of the estimate} $$\Big\|
\sum_{k=1}^{\infty} |dh_k| \Big\|_{1,\infty} \le \mathrm{c} \,
\|f\|_1.$$ Taking into account the form of $dh_k$ according to
Davis decomposition, we have
\begin{eqnarray*}
\Big\| \sum_{k=1}^\infty |dh_k| \Big\|_{1,\infty} & \le & 2 \Big\|
\sum_{k=1}^\infty \big| df_k \chi_{\{f_k^* \ge 2 f_{k-1}^*\}}
\big| \Big\|_{1,\infty} \\ & + & 2 \Big\| \sum_{k=1}^\infty \big|
\mathsf{E}_{k-1} \big( df_k \chi_{\{f_k^* \ge 2 f_{k-1}^*\}} \big)
\big| \Big\|_{1,\infty} = \mathrm{C} + \mathrm{D}.
\end{eqnarray*}
The estimate for the term $\mathrm{C}$ is very simple. Indeed, it
suffices to use the classical property of this part of Davis
decomposition. Namely, we have $f_k^* \ge 2 f_{k-1}^*$ if and only
if $f_k^* \le 2 (f_k^* - f_{k-1}^*)$. In particular, we deduce
$$\big| df_k \chi_{\{f_k^* \ge 2 f_{k-1}^*\}} \big| \le 2 f_k^*
\chi_{\{f_k^* \le 2 (f_k^* - f_{k-1}^*)\}} \le 4 (f_k^* -
f_{k-1}^*)$$ and conclude the following estimate $$\mathrm{C} \le
8 \sup_{\lambda > 0} \, \lambda \, \mu \Big\{ \sum_{k=1}^\infty
f_k^* - f_{k-1}^* > \lambda \Big\} = 8 \sup_{\lambda > 0} \,
\lambda \, \mu \big( f^* > \lambda \big) \le 8 \, \|f\|_1.$$ The
estimate for $\mathrm{D}$ is a little more complicated, we have
$$\mathrm{D} = \sup_{\lambda > 0} \, \lambda \, \mu \Big\{
\sum_{k=1}^\infty \big| \mathsf{E}_{k-1} \big( df_k \chi_{\{f_k^*
\ge 2 f_{k-1}^*\}} \big) \big|
> \lambda \Big\}.$$ Then we fix $\lambda > 0$ and apply Gundy's
decomposition \eqref{Eq-Gundy-Dec} to $f$ $$df_k = d\alpha_k + d
\beta_k + d \gamma_k.$$ By the quasi-triangle inequality we may
write
\begin{eqnarray*}
\lambda \, \mu \Big\{ \sum_{k=1}^\infty \big| \mathsf{E}_{k-1}
\big( df_k \chi_{\{f_k^* \ge 2 f_{k-1}^*\}} \big) \big| > \lambda
\Big\} \!\!\! & \le & \!\!\! \lambda \, \mu \Big\{
\sum_{k=1}^\infty \mathsf{E}_{k-1} \big( |d\alpha_k| \chi_{\{f_k^*
\ge 2 f_{k-1}^*\}} \big) > \frac{\lambda}{3} \Big\} \\ \!\!\! & +
& \!\!\! \lambda \, \mu \Big\{ \sum_{k=1}^\infty \mathsf{E}_{k-1}
\big( |d\beta_k| \chi_{\{f_k^* \ge 2 f_{k-1}^*\}} \big) >
\frac{\lambda}{3} \Big\}
\\ \!\!\! & + & \!\!\! \lambda \, \mu \Big\{ \sum_{k=1}^\infty \mathsf{E}_{k-1}
\big( |d\gamma_k| \chi_{\{f_k^* \ge 2 f_{k-1}^*\}} \big) >
\frac{\lambda}{3} \Big\}.
\end{eqnarray*}
We shall denote these terms by $\mathrm{D}_\alpha$,
$\mathrm{D}_\beta$ and $\mathrm{D}_\gamma$ respectively. The
estimates for $\alpha$ and $\gamma$ are straightforward. Indeed,
recalling that $d \alpha_k = df_k \chi_{\{f_{k-1}^*
> \lambda\}}$ and that the function $\chi_{\{f_{k-1}^* >
\lambda\}}$ is predictable, we conclude $$\mathrm{D}_\alpha =
\lambda \, \mu \Big\{ \sum_{k=1}^\infty \chi_{\{f_{k-1}^* >
\lambda\}} \mathsf{E}_{k-1} \big( |df_k| \chi_{\{f_k^* \ge 2
f_{k-1}^*\}} \big) > \frac{\lambda}{3} \Big\} \le \lambda \, \mu
\big( f^* > \lambda) \le \|f\|_1.$$ On the other hand, by
Chebychev's inequality we deduce
$$\mathrm{D}_\gamma \le 3 \Big\| \sum_{k=1}^\infty \mathsf{E}_{k-1}
\big( |d \gamma_k| \chi_{\{f_k^* \ge 2 f_{k-1}^*\}} \big) \Big\|_1
\le 3 \sum_{k=1}^\infty \|d \gamma_k\|_1 \le \mathrm{c} \,
\|f\|_1,$$ where the last inequality follows from property iii) in
Gundy's decomposition, see Section \ref{Section1}. It only remains
to estimate the term $\mathrm{D}_\beta$. To that aim, we first
recall the dual form of Doob's maximal inequality, since we shall
need it in the proof. Given $1 \le p < \infty$ and a sequence
$(\psi_m)_{m \ge 1}$ of positive functions in $L_p(\Omega)$, we
have
\begin{equation} \label{Eq-Dual-Doob}
\Big\| \summ_{m} \mathsf{E}_{m-1} (\phi_m) \Big\|_p \le \mathrm{c}
p \, \Big\| \summ_{m} \phi_m \Big\|_p.
\end{equation}
I learned this from Marius Junge, who obtained its noncommutative
analog in \cite{J1} as a key point to prove Doob's maximal
inequality for noncommutative martingales, see also Burkholder's
paper \cite{B} for more details. Applying Chebychev's inequality
for $p=3$ and the dual Doob inequality \eqref{Eq-Dual-Doob} for
the same index, we obtain
\begin{eqnarray*}
\mathrm{D}_\beta & \le & \mathrm{c} \lambda^{-2} \Big\|
\sum_{k=1}^\infty \mathsf{E}_{k-1} \big( |d\beta_k| \chi_{\{f_k^*
\ge 2 f_{k-1}^*\}} \big) \Big\|_3^3 \\ & \le & \mathrm{c}
\lambda^{-2} \Big\| \sum_{k=1}^\infty |d\beta_k| \chi_{\{f_k^* \ge
2 f_{k-1}^*\}} \Big\|_3^3 \\ & \le & \mathrm{c} \lambda^{-2}
\Big\| \Big( \sum_{k=1}^\infty |d\beta_k|^2 \Big)^{1/3} \Big(
\sum_{k=1}^\infty |d\beta_k|^{1/2} \chi_{\{f_k^* \ge 2
f_{k-1}^*\}} \Big)^{2/3} \Big\|_3^3 \\ & = & \mathrm{c}
\lambda^{-2} \Big\| \Big( \sum_{k=1}^\infty |d\beta_k|^2 \Big)
\Big( \sum_{k=1}^\infty |d\beta_k|^{1/2} \chi_{\{f_k^* \ge 2
f_{k-1}^*\}} \Big)^2 \Big\|_1 \\ & \le & \mathrm{c} \lambda^{-2}
\Big\| \sum_{k=1}^\infty |d\beta_k|^2 \Big\|_1 \, \Big\|
\sum_{k=1}^\infty |d\beta_k|^{1/2} \chi_{\{f_k^* \ge 2
f_{k-1}^*\}} \Big\|_\infty^2 \\ & \le & \mathrm{c} \lambda^{-1}
\Big\| \sum_{k=1}^\infty |d\beta_k|^{1/2} \chi_{\{f_k^* \ge 2
f_{k-1}^*\}} \Big\|_\infty^2 \|f\|_1,
\end{eqnarray*}
where the last inequality follows from the quadratic estimate in
property ii) of Gundy's decomposition, see Section \ref{Section1}.
Therefore, the only remaining estimate to conclude the proof of
Theorem A is the following
\begin{equation} \label{Eq-Linfty-Estimate}
\Big\| \sum_{k=1}^\infty |d\beta_k|^{1/2} \chi_{\{f_k^* \ge 2
f_{k-1}^*\}} \Big\|_\infty \le \mathrm{c} \sqrt{\lambda}.
\end{equation}
In order to prove inequality \eqref{Eq-Linfty-Estimate} we write
$d\beta_k$ as follows
\begin{eqnarray*}
d\beta_k & = & f_k \chi_{\{f_k^* \le \lambda\}} - f_{k-1}
\chi_{\{f_k^* \le \lambda\}}\\ & + & \mathsf{E}_{k-1} \big(
f_{k-1} \chi_{\{f_k^* \le \lambda\}} \big) - \mathsf{E}_{k-1}
\big( f_k \chi_{\{f_k^* \le \lambda\}} \big) = a_k-b_k+c_k-d_k.
\end{eqnarray*}
We clearly have
$$\max(a_k, b_k) \le f_k^* \chi_{\{f_k^* \le \lambda\}}.$$ On the
other hand, the following estimates hold
$$c_k = f_{k-1} \mathsf{E}_{k-1} (\chi_{\{f_k^* \le \lambda\}}) \le f_{k-1}
\mathsf{E}_{k-1} (\chi_{\{f_{k-1}^* \le \lambda\}}) \le f_{k-1}^*
\chi_{\{f_{k-1}^* \le \lambda\}},$$ $$d_k \le \mathsf{E}_{k-1}
\big( f_k \chi_{\{f_{k-1}^* \le \lambda\}} \big) = f_{k-1}
\chi_{\{f_{k-1}^* \le \lambda\}} \le f_{k-1}^* \chi_{\{f_{k-1}^*
\le \lambda\}}.$$ Thus we conclude $$|d \beta_k|^{1/2} \le
\sqrt{a_k + b_k} + \sqrt{c_k + d_k} \le \sqrt{2} \, \Big(
\sqrt{f_k^*} \chi_{\{f_k^* \le \lambda\}} + \sqrt{f_{k-1}^*}
\chi_{\{f_{k-1}^* \le \lambda\}} \Big).$$ This means that it
suffices to prove the following estimates
\begin{eqnarray}
\label{Eq-EstFinal1} \mathrm{E} & = & \Big\| \sum_{k=1}^\infty
\sqrt{f_k^*} \, \chi_{\{f_k^* \le \lambda\}} \, \chi_{\{f_k^* \ge
2 f_{k-1}^*\}} \Big\|_\infty \le \mathrm{c} \sqrt{\lambda}, \\
\label{Eq-EstFinal2} \mathrm{F} & = & \Big\| \sum_{k=1}^\infty
\sqrt{f_{k-1}^*} \, \chi_{\{f_{k-1}^* \le \lambda\}} \,
\chi_{\{f_k^* \ge 2 f_{k-1}^*\}} \Big\|_\infty \le \mathrm{c}
\sqrt{\lambda}.
\end{eqnarray}

\subsection{Step 3:} \textbf{Proof of the estimate} $$\max
\big( \mathrm{E}, \mathrm{F} \big) \le \frac{\sqrt{2}}{\sqrt{2}-1}
\, \sqrt{\lambda}.$$ Fix $w \in \Omega$ and set
$$\mathcal{S}_{\mathrm{E},n}(w) = \sum_{k=1}^n \phi_k(w) = \sum_{k=1}^n
\sqrt{f_k^*(w)} \, \chi_{\{f_k^* \le \lambda\}}(w) \,
\chi_{\{f_k^* \ge 2 f_{k-1}^*\}}(w).$$ We shall prove by induction
on $n$ that $\mathcal{S}_{\mathrm{E},n}(w) \le
\sqrt{2}/(\sqrt{2}-1) \, \sqrt{\lambda}$. Indeed, the assertion is
clear for $n=1$. Thus, let us assume that the assertion holds for
$n-1$ and let us estimate $\mathcal{S}_{\mathrm{E},n}(w)$. If
$\phi_n(w)=0$ we have $\mathcal{S}_{\mathrm{E},n}(w) =
\mathcal{S}_{\mathrm{E},n-1}(w)$ and there is nothing to prove. If
$\phi_n(w) \neq 0$, we must have $f_k^*(w) \le \lambda$ and
$$f_{n-1}^*(w) \le \frac12 f_n^*(w) \le \lambda/2.$$ Then, if
$\phi_{n-1}(w) \neq 0$ we know that $$f_{n-1}^*(w) \le \lambda/2
\quad \mbox{and} \quad f_{n-2}^*(w) \le \frac12 f_{n-1}^*(w) \le
\lambda/4.$$ On the other hand, if $\phi_{n-1}(w) = 0$ we may
ignore that term in the sum and we still have at our disposal that
$f_{n-2}^*(w) \le f_{n-1}^*(w) \le \lambda/2$, so that we can
argue in the same way for $\phi_{n-2}(w)$. Iterating the same
argument, it is not difficult to conclude that
$$\mathcal{S}_{\mathrm{E},n}(w) \le \sqrt{\lambda}
\sum_{k=0}^\infty \Big( \frac{1}{\sqrt{2}} \Big)^k =
\frac{\sqrt{2}}{\sqrt{2}-1} \, \sqrt{\lambda}.$$ Since this
argument works for any $w \in \Omega$, we have proved our claim
for $\mathrm{E}$. Arguing in a similar way, we obtain the same
bound for $\mathrm{F}$. Let us include the details for the sake of
completeness. Fix $w \in \Omega$ and set
$$\mathcal{S}_{\mathrm{F},n}(w) = \sum_{k=1}^n \psi_k(w) =
\sum_{k=1}^n \sqrt{f_{k-1}^*(w)} \, \chi_{\{f_{k-1}^* \le
\lambda\}}(w) \, \chi_{\{f_k^* \ge 2 f_{k-1}^*\}}(w).$$ Again, it
is clear that $\psi_1(w) \le \sqrt{\lambda}$ so that we may
proceed by induction on $n$ and assume that the inequality
$\mathcal{S}_{\mathrm{F},n-1}(w) \le \sqrt{2}/(\sqrt{2}-1) \,
\sqrt{\lambda}$ holds. If $\psi_n(w) = 0$ there is nothing to
prove while for $\psi_n(w) \neq 0$ we deduce that $f_{n-1}^*(w)
\le \lambda$. Going backwards, we seek for the next non-zero term
$\psi_j(w) \neq 0$. Such term implies $$f_{j-1}^*(w) \le \frac12
f_j^*(w) \le \frac12 f_{n-1}^*(w) \le \lambda/2.$$ Iterating one
more time we conclude $$\mathcal{S}_{\mathrm{F},n}(w) \le
\sqrt{\lambda} \sum_{k=0}^\infty \Big( \frac{1}{\sqrt{2}} \Big)^k
= \frac{\sqrt{2}}{\sqrt{2}-1} \, \sqrt{\lambda}.$$ This justifies
our claim $$\max \big( \mathrm{E}, \mathrm{F} \big) \le
\frac{\sqrt{2}}{\sqrt{2}-1} \, \sqrt{\lambda}$$ which implies
\eqref{Eq-EstFinal1} and \eqref{Eq-EstFinal2}. Therefore, the
proof of Theorem A is completed. \fin

\section{Proof of Theorem B}
\label{Section3}

Our aim in this section is proving Theorem B. However, instead of
giving the simplest proof available (see Paragraph
\ref{Paragraph3.4} for this), we present a more general proof
where the hypothesis of $\mathrm{k}$-regularity is only needed in
the very last step. This somehow supports our comment in the
Introduction on the validity of Theorem B for non-regular
martingales. Moreover, as we shall see below, our proof goes a
little further, see Corollary C below.

\vskip5pt

For now on and until the very end of the proof, we assume that $f$
is a bounded martingale in $L_1(\Omega)$ adapted to a
\emph{non-necessarily regular} filtration. Our martingale
decomposition $f = g + h$ is given one more time by Davis
decomposition
\begin{eqnarray*}
dg_k & = & df_k \chi_{\{f_k^* < 2 f_{k-1}^*\}} - \mathsf{E}_{k-1}
\big( df_k \chi_{\{f_k^* < 2 f_{k-1}^*\}} \big), \\ dh_k & = &
df_k \chi_{\{f_k^* \ge 2 f_{k-1}^*\}} - \mathsf{E}_{k-1} \big(
df_k \chi_{\{f_k^* \ge 2 f_{k-1}^*\}} \big).
\end{eqnarray*}
In particular, the weak type estimate
$$\Big\| \Big( \sum_{k=1}^\infty \mathsf{E}_{k-1}(|dg_k|^2)
\Big)^{1/2} \Big\|_{L_{1,\infty}(\Omega)} \le \mathrm{c} \,
\|f\|_1$$ holds with an absolute constant $\mathrm{c}$ by means of
Theorem A. In order to estimate the second term on the left of
\eqref{Eq-Weak-Burkholder2} we combine one more time Davis and
Gundy decompositions. More concretely, for fixed $\lambda
> 0$ and according to \eqref{Eq-Gundy-Dec} we have
\begin{eqnarray*}
dh_k & = & d\alpha_k \chi_{\{f_k^* \ge 2 f_{k-1}^*\}} -
\mathsf{E}_{k-1} \big( d\alpha_k \chi_{\{f_k^* \ge 2 f_{k-1}^*\}}
\big) \\ & + & d\beta_k \hskip0.5pt \chi_{\{f_k^* \ge 2
f_{k-1}^*\}} - \mathsf{E}_{k-1} \big( d\beta_k \chi_{\{f_k^* \ge 2
f_{k-1}^*\}} \big) \\ & + & d\gamma_k \, \chi_{\{f_k^* \ge 2
f_{k-1}^*\}} - \mathsf{E}_{k-1} \big( d\gamma_k \chi_{\{f_k^* \ge
2 f_{k-1}^*\}} \big) = dh_{\alpha k} + dh_{\beta k} + dh_{\gamma
k}.
\end{eqnarray*}
As in the proof of Theorem A, the terms associated to $\gamma$ are
the simplest ones
\begin{eqnarray*}
\Big\| \sum_{k=1}^\infty \delta_k \otimes d\gamma_k \,
\chi_{\{f_k^* \ge 2 f_{k-1}^*\}}
\Big\|_{L_{1,\infty}(\Omega_{\oplus \infty})} & \le &
\sum_{k=1}^\infty \|d\gamma_k\|_1 \le \mathrm{c} \, \|f\|_1, \\
\Big\| \sum_{k=1}^\infty \delta_k \otimes \mathsf{E}_{k-1} \big(
d\gamma_k \, \chi_{\{f_k^* \ge 2 f_{k-1}^*\}} \big)
\Big\|_{L_{1,\infty}(\Omega_{\oplus \infty})} & \le &
\sum_{k=1}^\infty \|d\gamma_k\|_1 \le \mathrm{c} \, \|f\|_1.
\end{eqnarray*}
Thus, it remains to estimate the terms associated to the
martingales $\alpha$ and $\beta$.

\subsection{Step 1:} \textbf{Proof of the estimate}
$$\lambda \sum_{k=1}^\infty \mu \Big\{ |dh_{\beta k}| > \lambda
\Big\} \le \mathrm{c} \, \|f\|_1.$$ This estimate is similar to
that of Theorem A. By Chebychev's inequality
\begin{eqnarray*}
\lambda \sum_{k=1}^\infty \mu \Big\{ |d\beta_k| \chi_{\{f_k^* \ge
2 f_{k-1}^*\}} > \lambda \Big\} & \le & \lambda^{-2}
\sum_{k=1}^\infty \big\| |d \beta_k| \, \chi_{\{f_k^* \ge 2
f_{k-1}^*\}} \big\|_3^3 \\ & = & \lambda^{-2} \Big\|
\sum_{k=1}^\infty |d \beta_k|^3 \, \chi_{\{f_k^* \ge 2
f_{k-1}^*\}} \Big\|_1 \\ & \le & \lambda^{-2} \Big\| \Big(
\sum_{k=1}^\infty |d \beta_k|^2 \Big) \, \sup_{k \ge 1} |d
\beta_k| \chi_{\{f_k^* \ge 2 f_{k-1}^*\}} \Big\|_1 \\ & \le & 4
\lambda^{-1} \Big\| \sum_{k=1}^\infty |d\beta_k|^2 \Big\|_1 =
\frac{4}{\lambda} \|\beta\|_2^2 \le \mathrm{c} \, \|f\|_1,
\end{eqnarray*}
where the last inequality follows from the quadratic estimate in
property ii) of Gundy's decomposition and the previous estimate
uses $|d\beta_k| \le 4 \lambda$ for all integer $k \ge 1$ since
$$|d \beta_k| = \Big| df_k \chi_{\{f_k^* \le \lambda\}} -
\mathsf{E}_{k-1} \big( df_k \chi_{\{f_k^* \le \lambda\}} \big)
\Big| \le 2 f_k^* \chi_{\{f_k^* \le \lambda\}} + 2
\mathsf{E}_{k-1} \big( f_k^* \chi_{\{f_k^* \le \lambda\}} \big)
\le 4 \lambda.$$ In a similar way, since the $\mathsf{E}_{k-1}$'s
are contractive in $L_3(\Omega)$, we find $$\lambda
\sum_{k=1}^\infty \mu \Big\{ \mathsf{E}_{k-1} \big( |d\beta_k|
\chi_{\{f_k^* \ge 2 f_{k-1}^*\}} \big) > \lambda \Big\} \le
\mathrm{c} \, \|f\|_1.$$ The assertion follows from the estimates
above and the quasi-triangle inequality.

\subsection{Step 2:} \textbf{Proof of the estimate} $$\lambda
\sum_{k=1}^\infty \mu \Big\{ |d\alpha_k| \chi_{\{f_k^* \ge 2
f_{k-1}^*\}} > \lambda \Big\} \le \mathrm{c} \, \|f\|_1.$$ Let us
fix $\lambda > 0$. Since $d\alpha_k = df_k \chi_{\{f_{k-1}^* >
\lambda\}}$, we have
$$\lambda \sum_{k=1}^\infty \mu \Big\{ |d\alpha_k| \chi_{\{f_k^*
\ge 2 f_{k-1}^*\}} > \lambda \Big\} \le \lambda \sum_{k=1}^\infty
\mu \Big\{ f_k^* \ge 2 f_{k-1}^* > 2 \lambda \Big\}.$$ Then it is
clear from Doob's maximal inequality that the set of points $w \in
\Omega$ such that $f_k^*(w) \ge 2 f_{k-1}^*(w) > 2 \lambda$ for
infinitely many $k$'s has zero $\mu$-measure. Therefore, denoting
by $\mathcal{P}_f(\N)$ the (numerable) set of subsets of $\N$ with
finite cardinal, we may decompose the sum
$$\lambda \sum_{k=1}^\infty \mu \Big\{ f_k^* \ge 2 f_{k-1}^* > 2
\lambda \Big\}$$ as follows
$$\lambda \sum_{s=1}^\infty s \sum_{\begin{subarray}{c}
\mathcal{A} \in \mathcal{P}_f(\N) \\ \# \mathcal{A}=s
\end{subarray}} \mu \Big\{ f_k^* \ge 2 f_{k-1}^* > 2 \lambda
\Leftrightarrow k \in \mathcal{A} \Big\} = \lambda
\sum_{s=1}^\infty s \sum_{\begin{subarray}{c} \mathcal{A} \in
\mathcal{P}_f(\N) \\ \# \mathcal{A}=s
\end{subarray}} \mu (\Omega_{\mathcal{A},\lambda}).$$ It is clear
from the definition that $\Omega_{\mathcal{A}_1,\lambda} \cap
\Omega_{\mathcal{A}_2,\lambda} = \emptyset$ with $\mathcal{A}_1
\neq \mathcal{A}_2$. Moreover, we also recall that each
$\Omega_{\mathcal{A},\lambda}$ with $\# \mathcal{A} = s$ is
contained in the set where $f^* > 2^s \lambda$ since the condition
$f_k^* \ge 2 f_{k-1}^*$ is satisfied $s$ times and $f_{k_0-1}^* >
\lambda$ for the smallest integer $k_0 \in \mathcal{A}$. Hence, we
have $$\bigsqcup_{\begin{subarray}{c} \mathcal{A} \in
\mathcal{P}_f(\N)
\\ \# \mathcal{A}=s
\end{subarray}} \Omega_{\mathcal{A},\lambda} \subset \Big\{ f^* >
2^s \lambda \Big\}.$$ By means of Doob's maximal inequality we
conclude $$\lambda \sum_{s=1}^\infty s \sum_{\begin{subarray}{c}
\mathcal{A} \in \mathcal{P}_f(\N) \\ \# \mathcal{A}=s
\end{subarray}} \mu (\Omega_{\mathcal{A},\lambda}) \le \lambda
\sum_{s=1}^\infty s \mu \big( f^* > 2^s \lambda \big) \le \Big(
\sum_{s=1}^\infty s/2^s \Big) \, \|f\|_1 = 2 \, \|f\|_1.$$

\subsection{Step 3:} \textbf{Proof of the remaining estimate for
$\mathrm{k}$-regular martingales}
\begin{equation} \label{Eq-Key-Estimate}
\lambda \sum_{k=1}^\infty \mu \Big\{ \mathsf{E}_{k-1} \big(
|d\alpha_k| \chi_{\{f_k^* \ge 2 f_{k-1}^*\}} \big) > \lambda
\Big\} \le \mathrm{c} \mathrm{k} \, \|f\|_1.
\end{equation}
This is the only estimate in the proof of Theorem B where we use
the $\mathrm{k}$-regularity of the filtration $\mathsf{A}_1,
\mathsf{A}_2, \ldots$ Let us introduce some notation. Let us
consider a parameter $\mathrm{k} > 1$. We shall say that the
filtration $\mathsf{A}_1, \mathsf{A}_2, \ldots$ is
\emph{$\mathrm{k}$-homogeneous} if for every $n \ge 1$ and any
measurable set $\mathrm{A} \in \mathsf{A}_n$ we have $$\mu \Big\{
\mbox{supp} \, \mathsf{E}_{n-1} (\chi_{\mathrm{A}}) \Big\} \le
\mathrm{k} \, \mu(\mathrm{A}).$$ We are not aware if this notion
appears somewhere in the literature. In the following result, we
study the relation between $\mathrm{k}_1$-regularity and
$\mathrm{k}_2$-homogeneity. Again, we ignore whether or not this
result is already known.

\begin{lemma} \label{Lemma-Reg-Hom}
Given $1 < \mathrm{k}_1 < \mathrm{k}_2 < \infty$, we have
$$\mbox{$\mathrm{k}_1$-regularity} \ \Rightarrow \
\mbox{$\mathrm{k}_1$-homogeneity} \ \Rightarrow \
\mbox{$\mathrm{k}_2$-regularity}.$$
\end{lemma}

\dem Let us begin with the first implication. Given a positive
integer $n \ge 1$ and a measurable set $\mathrm{A} \in
\mathsf{A}_n$, we known from the assumption on
$\mathrm{k}_1$-regularity that the inequality below holds
\begin{equation} \label{Eq-k1-Regular}
\chi_\mathrm{A} \le \mathrm{k}_1 \, \mathsf{E}_{n-1}
(\chi_\mathrm{A}).
\end{equation}
Let us assume that $$\mu \Big\{ \mbox{supp} \, \mathsf{E}_{n-1}
(\chi_\mathrm{A}) \Big\} > \mathrm{k}_1 \, \mu(\mathrm{A}).$$ By
the main property of the conditional expectation, we have
\begin{eqnarray*}
\mathrm{k}_1 \, \mu(\mathrm{A}) & = & \mathrm{k}_1 \, \int_\Omega
\chi_\mathrm{A} \, d \mu \\ & = & \mathrm{k}_1 \, \int_\Omega
\mathsf{E}_{n-1} (\chi_\mathrm{A}) \, d \mu \\ & \ge &
\mathrm{k}_1 \, \mu \Big\{ \mbox{supp} \,
\mathsf{E}_{n-1}(\chi_\mathrm{A}) \Big\} \, \inf \Big\{
\mathsf{E}_{n-1}(\chi_\mathrm{A})(w) \, \big| \, w \in
\mathrm{supp} \, \mathsf{E}_{n-1}(\chi_\mathrm{A}) \Big\}.
\end{eqnarray*}
In particular, we deduce that $$\inf \Big\{
\mathsf{E}_{n-1}(\chi_\mathrm{A})(w) \, \big| \, w \in
\mathrm{supp} \, \mathsf{E}_{n-1}(\chi_\mathrm{A}) \Big\} <
1/\mathrm{k}_1$$ so that the following set in $\mathsf{A}_{n-1}$
$$\mathrm{B} = \Big\{ 0 < \mathsf{E}_{n-1} (\chi_\mathrm{A}) <
1/\mathrm{k}_1 \Big\}$$ has positive $\mu$-measure. If we show
that $\mu \big( \mathrm{A} \cap \mathrm{B} \big) > 0$ we will
conclude since $\mathrm{A} \cap \mathrm{B}$ is a subset of
$\mathrm{A}$ where \eqref{Eq-k1-Regular} fails. Let us assume that
$\mu \big( \mathrm{A} \cap \mathrm{B} \big) = 0$. In that case we
have $\chi_\mathrm{A} \chi_\mathrm{B} = 0$ $\mu$-a.e. and since
$\mathrm{B} \in \mathsf{A}_{n-1}$ we conclude $$\mathsf{E}_{n-1}
(\chi_\mathrm{A}) \chi_{\mathrm{B}} = \mathsf{E}_{n-1}
(\chi_\mathrm{A} \chi_{\mathrm{B}}) = 0 \quad \mbox{$\mu$-a.e.}$$
However, this contradicts the definition of $\mathrm{B}$ and the
proof of the first implication is completed. For the second
implication, we have to show that the $\mathrm{A}_n$-measurable
set $\mathrm{A} = \big\{ f_n
> \mathrm{k}_2 f_{n-1} \big\}$ has zero $\mu$-measure. Let us
assume that $\mu(\mathrm{A}) > 0$. Then we consider a
$\mathsf{A}_n$-measurable function $g$ of the form $g = \sum_j
\alpha_j \chi_{\mathrm{B}_j}$, where the sum might have infinitely
many terms, the $\mathrm{B}_j$'s are in $\mathsf{A}_n$ and such
that $$g \le f_n \quad \mbox{and} \quad \|g - f_n\|_\infty <
\varepsilon.$$ If $w \in \mathrm{A}$, we have $$g(w) > f_n(w) -
\varepsilon > \mathrm{k}_2 f_{n-1}(w) - \varepsilon \ge
\mathrm{k}_2 \mathsf{E}_{n-1}(g)(w) - \varepsilon.$$ Therefore it
turns out that $$\mathrm{A} \subset \Big\{ \summ_j \alpha_j
\chi_{\mathrm{B}_j} > \mathrm{k}_2 \summ_j \alpha_j
\mathsf{E}_{n-1} ( \chi_{\mathrm{B}_j}) - \varepsilon \Big\}
\subset \bigcup_{j \ge 1} \Big\{ \chi_{\mathrm{B}_j} >
\mathrm{k}_2 \mathsf{E}_{n-1}(\chi_{\mathrm{B}_j}) - \varepsilon
\Big\}.$$ Recalling that $\mu(\mathrm{A}) > 0$, there must exist
some $j_0 \ge 1$ such that $$\mu \Big\{ \chi_{\mathrm{B}_{j_0}} >
\mathrm{k}_2 \mathsf{E}_{n-1}(\chi_{\mathrm{B}_{j_0}}) -
\varepsilon \Big\} > 0.$$ We shall denote this set in
$\mathsf{A}_n$ by $\mathrm{D}_{j_0}$. Given $w \in
\mathrm{D}_{j_0}$, we clearly have
\begin{equation} \label{Eq-Dj0}
\mathrm{k}_2 \mathsf{E}_{n-1}(\chi_{\mathrm{D}_{j_0}})(w) < 1 +
\varepsilon.
\end{equation}
On the other hand, by $\mathrm{k}_1$-homogeneity
\begin{eqnarray*}
\mu \Big\{ \mbox{supp} \, \mathrm{E}_{n-1}
(\chi_{\mathrm{D}_{j_0}}) \Big\} & \le & \mathrm{k}_1 \, \mu
(\mathrm{D}_{j_0}) \\ & = & \mathrm{k}_1 \int_\Omega
\mathsf{E}_{n-1}(\chi_{\mathrm{D}_{j_0}}) \, d \mu
\\ & \le & \mathrm{k}_1 \, \mu \Big\{ \mbox{supp} \, \mathrm{E}_{n-1}
(\chi_{\mathrm{D}_{j_0}}) \Big\} \, \big\| \mathsf{E}_{n-1}
(\chi_{\mathrm{D}_{j_0}}) \big\|_\infty.
\end{eqnarray*}
This means that $\big\| \mathsf{E}_{n-1} (\chi_{\mathrm{D}_{j_0}})
\big\|_\infty \ge 1 / \mathrm{k}_1$, hence
$$\mu(\mathrm{R}_\delta) > 0 \quad \mbox{with} \quad
\mathrm{R}_\delta = \Big\{ \mathsf{E}_{n-1}
(\chi_{\mathrm{D}_{j_0}}) > 1/\mathrm{k}_1 - \delta \Big\} \quad
\mbox{for every $\delta > 0$}.$$ Now we observe that for any $w
\in \mathrm{D}_{j_0} \cap \mathrm{R}_\delta$ we obtain
$$\frac{\mathrm{k}_2}{\mathrm{k}_1} - \delta \mathrm{k}_2 <
\mathrm{k}_2 \mathsf{E}_{n-1}(\chi_{\mathrm{D}_{j_0}})(w) < 1 +
\varepsilon.$$ Therefore, since we are assuming that $\mathrm{k}_1
< \mathrm{k}_2$, we may take $\delta$ and $\varepsilon$ small
enough so that the relation above provides the desired
contradiction as far as we show that $\mu(\mathrm{D}_{j_0} \cap
\mathrm{R}_\delta) > 0$. However, this follows as in the first
part of the proof. Indeed, if not we would have
$\chi_{\mathrm{D}_{j_0}} \chi_{\mathrm{R}_\delta} = 0$ $\mu$-a.e.
Thus, since $\mathrm{R}_\delta \in \mathsf{A}_{n-1}$ we should
conclude that $\mathsf{E}_{n-1}(\chi_{\mathrm{D}_{j_0}})
\chi_{\mathrm{R}_\delta} = 0$ $\mu$-a.e., which contradicts the
definition of $\mathrm{R}_{\delta}$. \fin

Now it is straightforward to finish the proof of Theorem B.
Namely, it remains to prove inequality \eqref{Eq-Key-Estimate}.
However, since we assume that $f$ is a positive martingale, we
have by $\mathrm{k}$-regularity that $df_k \le (\mathrm{k}-1)
f_{k-1}$. Therefore, assuming $\mathrm{k}$-regularity and
according to the first half of Lemma \ref{Lemma-Reg-Hom} we obtain
\begin{eqnarray*}
\lambda \sum_{k=1}^\infty \mu \Big\{ \mathsf{E}_{k-1} \big(
|d\alpha_k| \chi_{\{f_k^* \ge 2 f_{k-1}^*\}} \big)
> \lambda \Big\} \!\!\! & \le & \!\!\!
\lambda \sum_{k=1}^\infty \mu \Big\{ \mathsf{E}_{k-1} \big( \,
\chi_{\{f_k^* \ge 2 f_{k-1}^* > 2 \lambda\}} \big) > 0 \Big\}
\\ \!\!\! & \le & \!\!\! \mathrm{k} \lambda \sum_{k=1}^\infty \mu
\Big\{ f_k^* \ge 2 f_{k-1}^*
> 2 \lambda \Big\} \le 2 \mathrm{k} \, \|f\|_1,
\end{eqnarray*}
where the last inequality follows by Step 2. Theorem B is proved.
\fin

\subsection{Further remarks}
\label{Paragraph3.4}

We conclude this section by analyzing Theorem B in some detail. If
we do not care about its validity for non-regular martingales, a
much simpler proof is available. The idea is that regular
martingales are also characterized by being \emph{previsible}
martingales, see e.g. Proposition 2.19 in \cite{W}. This can be
used to show that the only relevant part in
\eqref{Eq-Weak-Burkholder2} for regular martingales is the term
associated to the conditional square function.

\vskip5pt

\noindent \textbf{A simpler proof of Theorem B.} Let $f = (f_1,
f_2, \ldots)$ be a bounded positive martingale in $L_1(\Omega)$
adapted to a $\mathrm{k}$-regular filtration. Then we may consider
the Davis type decomposition $f=g+h$ with martingale differences
given by
\begin{eqnarray*}
dg_k & = & df_k \chi_{\{f_k \le \mathrm{k} f_{k-1}\}} -
\mathsf{E}_{k-1} \big( df_k \chi_{\{f_k \le \mathrm{k} f_{k-1}\}}
\big), \\ dh_k & = & df_k \chi_{\{f_k > \mathrm{k} f_{k-1}\}} -
\mathsf{E}_{k-1} \big( df_k \chi_{\{f_k > \mathrm{k} f_{k-1}\}}
\big).
\end{eqnarray*}
According to the $\mathrm{k}$-regularity, the only non-zero
martingale difference in $h$ is $dh_1$. Thus, the second term in
\eqref{Eq-Weak-Burkholder2} is trivially controlled by $4\|f\|_1$.
On the other hand, the first term in \eqref{Eq-Weak-Burkholder2}
can be estimated as in Step 1 of the proof of Theorem A, with the
only difference that we obtain the constant $\mathrm{ck}$ instead
of $\mathrm{c}$. \fin

\begin{remark}
\emph{The proof given above shows that in the $\mathrm{k}$-regular
case it suffices to consider the conditional term and ignore the
diagonal one. As it was justified in Remark 8.3 of \cite{BG}, this
is only possible under the assumption of regularity.}
\end{remark}

As we already mentioned in the Introduction, martingale
inequalities where a martingale decomposition is involved arise
very naturally in the noncommutative setting, mainly due to the
row/column nature of the corresponding martingale Hardy spaces.
Among many other papers, we refer the reader to \cite{LuP,PX1,X}
for some illustrations of this phenomenon. However, not requiring
the decompositions to be martingale decompositions, it is
sometimes simpler to obtain the corresponding inequality
\cite{R2,R4}. Our first proof of Theorem B goes a little further
and produces the following result in this line.

\begin{corC}
Let $f = (f_1, f_2, \ldots)$ be a bounded martingale in
$L_1(\Omega)$. Then, we can decompose each $f_n$ as a sum $f_n =
g_n + h_n$ of two functions $($non-necessarily martingales$)$
adapted to the same filtration and satisfying the following
inequality with an absolute constant $\mathrm{c}$ $$\Big\| \Big(
\sum_{k=1}^\infty \mathsf{E}_{k-1}(|dg_k|^2) \Big)^{1/2}
\Big\|_{L_{1,\infty}(\Omega)} + \Big\| \sum_{k=1}^\infty \delta_k
\otimes dh_k \Big\|_{L_{1,\infty}(\Omega_{\oplus \infty})} \le
\mathrm{c} \, \|f\|_1.$$
\end{corC}

\dem Let us consider the decomposition
\begin{eqnarray*}
dg_k & = & df_k \chi_{\{f_k^* < 2 f_{k-1}^*\}}, \\ dh_k & = & df_k
\chi_{\{f_k^* \ge 2 f_{k-1}^*\}}.
\end{eqnarray*}
The result follows from our first proof of Theorem B since
\eqref{Eq-Key-Estimate} is not needed. \fin

Namely, we can drop the $\mathrm{k}$-regularity assumption as far
as we do not require to have a decomposition of $f$ into two
martingales. We have proved Theorem B under the assumption of
$\mathrm{k}$-regularity or avoiding martingale decompositions as
in Corollary C. However, it is still open to decide whether
Theorem B holds for arbitrary martingales. Let us state this
problem for the interested reader.

\begin{problem}
Let $f = (f_1, f_2, \ldots)$ be a bounded martingale in
$L_1(\Omega)$. Is there a decomposition of $f$ as a sum $f= g+h$
of two martingales adapted to the same filtration and satisfying
the following inequality with an absolute constant $\mathrm{c}$?
$$\Big\| \Big( \sum_{k=1}^\infty \mathsf{E}_{k-1}(|dg_k|^2)
\Big)^{1/2} \Big\|_{L_{1,\infty}(\Omega)} + \Big\|
\sum_{k=1}^\infty \delta_k \otimes dh_k
\Big\|_{L_{1,\infty}(\Omega_{\oplus \infty})} \le \mathrm{c} \,
\|f\|_1.$$
\end{problem}

\section{Applications and Comments}
\label{Section4}

In this section, we obtain some applications of Theorems A and B
and analyze their relation to the the theory of noncommutative
martingales. Let us begin by studying the implications of Theorem
A in Davis decomposition.

\subsection{On the classical Davis decomposition}

Besides its clear relation with Burkholder's martingale
inequality, Theorem A can also be understood as a weak type
estimate which generalizes the known properties of Davis
decomposition. In other words, recalling the second estimate in
Davis decomposition \eqref{Eq-Properties-Davis}
$$\Big\| \sum_{k=1}^\infty |dh_k| \Big\|_p \le \mathrm{c} p \,
\|f^*\|_p \le \mathrm{c} \frac{p^2}{p-1} \, \|f\|_p \quad
\mbox{for} \quad 1 < p < \infty,$$ the inequality
$$\Big\| \sum_{k=1}^\infty |dh_k| \Big\|_{1,\infty} \le \mathrm{c}
\, \|f\|_1$$ can be regarded as the associated weak type
inequality for $p=1$. This inequality was justified in Steps 2 and
3 of our proof of Theorem A. On the other hand, in a less explicit
way, the estimate for the first term in \eqref{Eq-Weak-Burkholder}
can also be understood as an extension of the known estimates for
$g$ in Davis decomposition.

\subsection{Real interpolation}

One of the first applications of our results that comes to mind is
reproving Burkholder's inequality by real interpolation and
duality. This is even possible starting from Corollary C, as
showed in \cite{R4}. This alternative proof has given rise to the
optimal constants for the noncommutative Burkholder inequality
\cite{JX}. It can be easily checked that the constants obtained in
\cite{R4} are still optimal in the commutative case as $p \to 1$
but not as $p \to \infty$, as it follows from Hitczenko's results
\cite{H}.

\vskip5pt

We should also point out that, in contrast with the previous
paragraph (where Theorem A came into scene), Theorem B is the
right result to obtain Burkholder's inequality to apply real
interpolation. Indeed, if we want to obtain \eqref{Eq-Burk-p<2}
from Theorem B, we just need the well-known isomorphism
$$\big[ L_{1,\infty}(\Omega_{\oplus \infty}), L_2(\Omega_{\oplus
\infty}) \big]_{\theta,p} \sim_{\mathrm{c}_p} L_p(\Omega_{\oplus
\infty}) \quad \mbox{with} \quad \mathrm{c}_p \sim p / p-1.$$ On
the contrary, Theorem A would require the interpolation
isomorphism
\begin{equation} \label{Eq-RealInt}
\big[ L_{1,\infty}(\Omega;\ell_1), L_2(\Omega;\ell_2)
\big]_{\theta,p} \sim_{\mathrm{c}_p} L_p(\Omega;\ell_p).
\end{equation}
The question whether \eqref{Eq-RealInt} holds or not is
reminiscent of other natural questions which people have asked
over the years with regard to extending the phenomenon of the
Marcinkiewicz interpolation theorem to more general situations:
bilinear operators, changes of weights or in our case,
vector-valued function spaces. Most conjectures about extending
Marcinkiewicz theorem in these various ways have turned out to be
untrue. However, after consulting several experts in the field the
question seemed quite unclear in this particular case. Finally,
Michael Cwikel showed to me \cite{Cw} that the isomorphism stated
above is also \emph{false}. Let us include Cwikel's argument for
completeness. Let us define $$\mathcal{A}_0 = L_{1,\infty}(\Omega;
\ell_1) \quad \mbox{and} \quad \mathcal{A}_1 = L_2(\Omega;
\ell_2).$$ It is clear that $\mathcal{A}_0$ contains
$L_1(\Omega_{\oplus_\infty})$, so that the inclusion below holds
$$L_p(\Omega;\ell_p) = L_p(\Omega_{\oplus_\infty}) \subset
[\mathcal{A}_0, \mathcal{A}_1]_{\theta,p}.$$ However, the reverse
inclusion fails. Following Cwikel's argument we show that the even
smaller space $\mathcal{A}_0 \cap \mathcal{A}_1$ is not contained
in $L_p(\Omega; \ell_p)$. To that aim we take $\Omega$ to be the
unit interval equipped with the Lebesgue measure. We will
represent elements of $\mathcal{A}_0 \cap \mathcal{A}_1$ as step
functions $f: [0,1] \times \R_+ \to \C$ of two variables $(w,z)$.
Of course, we need to assume that for each constant $w$, the
function $f(w, \cdot)$ is a constant function of $z$ on the
interval $(n-1,n]$ for each positive integer $n$. Let us consider
the following sets in the plane
$$\mathcal{S} = \bigcup_{k \ge 1} \mathcal{S}_k \quad \mbox{with}
\quad \mathcal{S}_k = \Big\{ 0 \le w \le 1/k \ \ \mbox{and} \ \
k-1 < z \le k \Big\}.$$ The set $\mathcal{S}$ is a
\emph{discretized} version of the region where $0 < z \le 1/w$.
Let $f_1$ be the characteristic function of $\mathcal{S}$. We
claim that, under the representation considered above, $f_1$ is an
element of $\mathcal{A}_0$. Indeed, we have
$$\|f_1\|_{\mathcal{A}_0} = \sup_{\lambda > 0} \lambda \, \mu \Big\{
\summ_k |f_1(w,k)| > \lambda \Big\} = 1 < \infty.$$ Now, let
$\alpha > 0$ be a positive number and let $$f_2(w,k) =
\frac{\chi_{\mathcal{S}_k}(w)}{(1 + \log k)^\alpha} \quad
\mbox{for} \quad k \ge 1.$$ Since $f_2 \le f_1$, it follows that
$f_2 \in \mathcal{A}_0$. Moreover, we have
$$\|f_2\|_{L_p(\Omega;\ell_p)} = \Big( \summ_k \Big\|
\frac{\chi_{\mathcal{S}_k}(w)}{(1 + \log k)^\alpha} \Big\|_p^p
\Big)^{1/p} = \Big( \summ_k \frac{1}{k (1 + \log k)^{p\alpha}}
\Big)^{1/p},$$ so that $f_2 \in L_p(\Omega, \ell_p)$ if and only
if $p\alpha > 1$. Therefore, given any $1 < p < 2$ we may choose
$\alpha$ so that $p \alpha < 1 < 2 \alpha$. In this case we obtain
that $$f_2 \in \mathcal{A}_0 \cap \mathcal{A}_1 \setminus
L_p(\Omega; \ell_p).$$

\subsection{A noncommutative Davis type decomposition}

In this paragraph we present a useful way to generalize the
martingale differences $dg_k$ and $dh_k$ in Davis martingale
decomposition to the noncommutative setting. We shall assume
certain familiarity with the theory of noncommutative martingale
inequalities. We begin with the trivial identity
\begin{eqnarray*}
\Big\{ f_k^* < 2 f_{k-1}^* \Big\} & = & \bigcup_{\lambda > 0}
\Big\{ \frac{\lambda}{2} < f_{k-1}^* \le \lambda \Big\} \cap
\Big\{ \frac{\lambda}{2} < f_k^* \le \lambda \Big\} \\ & = &
\bigcup_{\lambda > 0} \Big( \Big\{ f_{k-1}^* \le \lambda \Big\}
\setminus \Big\{ f_{k-1}^* \le \frac{\lambda}{2} \Big\} \Big) \cap
\Big( \Big\{ f_k^* \le \lambda \Big\} \setminus \Big\{ f_k^* \le
\frac{\lambda}{2}\Big\} \Big).
\end{eqnarray*}
One important lack in noncommutative martingales is the absence of
stopping times or maximal functions. In the case of maximal
functions, there exist two natural substitutes. Roughly speaking,
we use a construction due to Cuculescu \cite{C} when dealing with
weak type inequalities, while for strong inequalities the right
notion was formulated by Junge in \cite{J1}. We shall use here
Cuculescu's construction. In other words, for any $\lambda > 0$ we
can construct a sequence of projections $q_1(\lambda),
q_2(\lambda), \ldots$ which play the role of the sets $\{f_k^* \le
\lambda\}$ for $k=1,2,\ldots$ Namely, the sequence of
$q_k(\lambda)$'s is adapted for all $\lambda > 0$ and satisfies an
analogue of the weak type $(1,1)$ Doob's maximal inequality, see
\cite{C,PR,R} for more details. There is however one natural
property which is not satisfied by Cuculescu's projections. In
contrast with the classical case, where we have $\{f_k^* \le
\lambda_1\} \subset \{f_k^* \le \lambda_2\}$ whenever $\lambda_1 <
\lambda_2$, it is no longer true that $q_k(\lambda_1)$ is a
subprojection of $q_k(\lambda_2)$. This is solved by defining the
projections $$\pi_k(\lambda) = \bigwedge_{\xi \ge \lambda}
q_k(\xi) \sim \bigcap_{\xi \ge \lambda} \Big\{ f_k^* \le \xi
\Big\} = \Big\{ f_k^* \le \lambda \Big\}.$$

\vskip5pt

However, for some technical reasons like commuting properties of
the resulting projections, it is better to consider countable
families of $q_k$'s. More concretely, we consider dyadic
$\lambda$'s of the form $\lambda = 2^j$ and define
$$\pi_k(\lambda) = \bigwedge_{s=0}^\infty q_k(2^{j+s}) \sim \Big\{
f_k^* \le 2^j \Big\}.$$ This gives rise to the following
\emph{approximation}
\begin{eqnarray*}
\Big\{ f_k^* < 2 f_{k-1}^* \Big\} & \sim & \bigvee_{j=0}^\infty
\Big( \pi_{k-1}(2^j) - \pi_{k-1}(2^{j-1}) \Big) \wedge \Big(
\pi_k(2^j) - \pi_k(2^{j-1}) \Big) \\ & = & \sum_{j=0}^\infty \Big(
\pi_{k-1}(2^j) - \pi_{k-1}(2^{j-1}) \Big) \wedge \Big( \pi_k(2^j)
- \pi_k(2^{j-1}) \Big),
\end{eqnarray*}
where the last identity follows by pairwise orthogonality. Let us
observe in passing that the information lost in the process of
taking only dyadic $\lambda$'s is not relevant since our proofs of
Theorems A and B still hold when replacing the sets $\{f_k^* < 2
f_{k-1}^*\}$ by the smaller ones $$\bigcup_{j=0}^\infty \Big\{
2^{j-1} < f_{k-1}^* \le f_k^* \le 2^j \Big\}.$$

We have already rewritten the Davis sets $\{f_k^* < 2 f_{k-1}^*\}$
in a way that works in the noncommutative setting. Of course, this
gives rise to a Davis type decomposition for noncommutative
martingales. Although not mentioned by Randrianantoanina, this
might be a good motivation for the decomposition used in
\cite{R4}. Indeed, although the decomposition in \cite{R4} is not
made of martingale differences, the same idea is used there in a
very indirect way. More concretely, the projections $$p_{j,k} =
\bigwedge_{s=j}^\infty q_k(2^s) - \bigwedge_{s=j-1}^\infty
q_k(2^s) = \pi_k(2^j) - \pi_k(2^{j-1})$$ as well as $p_{j,k-1}
p_{j,k}$ are important key tools in \cite{R4}. Nevertheless, it is
also worthy of mention that the decomposition used there is not
exactly the \emph{translation} (via the transformations described
in this paragraph) of this paper. Namely, an extra
\emph{nonsymmetric} row/column partition is needed in the
noncommutative case.

\vskip5pt

Finally, we observe that the \emph{noncommutative Davis
decomposition} presented here is constructed with the aim to
interact with Cuculescu's construction. In other words, according
to the \emph{philosophy} mentioned above, this decomposition
should be the right one when dealing with weak type inequalities.
For strong inequalities, we should work with Junge's approach
\cite{J1}, but we still do not know how to obtain the right
\emph{maximal operators} and Davis decomposition is still unclear.

\subsection{Related results}

This paper is strongly motivated by problems and results from
noncommutative probability. We conclude by giving some references
(not included in the Introduction) related to this paper. The
optimal growth of the relevant constants in several noncommutative
martingale inequalities can be found in \cite{JX2}. Gundy's
decomposition of noncommutative martingales was obtained in the
recent paper \cite{PR}. There also exists free analogs of
generalized Khintchine and Rosenthal inequalities, see \cite{JPX}
and the references therein. We finally refer to \cite{J2,JP} for
applications of these results in the theory of noncommutative
$L_p$ spaces.

\renewcommand{\theequation}{\arabic{equation}}

\bibliographystyle{amsplain}

\begin{thebibliography}{10}
\bibitem {B} D.L. Burkholder, \emph{Distribution function
inequalities for martingales}. Ann. Probab. \textbf{1} (1973),
19-42.
\bibitem {BG} D.L. Burkholder and R.F. Gundy, \emph{Extrapolation
and interpolation of quasi-linear operators on martingales}. Acta
Math. \textbf{124} (1970), 249-304.
\bibitem {C} I. Cuculescu, \emph{Martingales on von Neumann
algebras}. J. Multivariate Anal. \textbf{1} (1971), 17-27.
\bibitem {Cw} M. Cwikel, Personal communication.
\bibitem {D} B. Davis, \emph{On the integrability of the martingale
square function}. Israel J. Math. \textbf{8} (1970), 187-190.
\bibitem {GA} A.M. Garsia, Martingale inequalities: {S}eminar
notes on recent progress. W. A. Benjamin, Inc., Reading,
Mass.-London-Amsterdam, 1973, Mathematics Lecture Notes Series.
\bibitem {G} R.F. Gundy, \emph{A decomposition for $L^1$-bounded
martingales}. Ann. Math. Statist. \textbf{39} (1968), 134-138.
\bibitem {H} P. Hitczenko, \emph{Best constants in martingale
version of Rosenthal's inequality}. Ann. Probab. \textbf{18}
(1990), 1656-1668.
\bibitem {JSZ} W.B. Johnson, G. Schechtman and J. Zinn, \emph{Best
constants in moment inequalities for linear combinations of
independent and exchangeable random variables}. Ann. Probab.
\textbf{13} (1985), 234-253.
\bibitem {J1} M. Junge, \emph{Doob's inequality for non-commutative
martingales}. J. reine angew. Math. \textbf{549} (2002), 149-190.
\bibitem {J2} M. Junge, \emph{Embedding of the operator space
$\mathrm{OH}$ and the logarithmic \lq little Grothendieck
inequality\rq}. Invent. Math. \textbf{161} (2005), 225-286.
\bibitem {JP} M. Junge and J. Parcet, \emph{Intersections and
interpolation in noncommutative $L_p$ spaces}. Preprint 2005.
\bibitem {JPX} M. Junge, J. Parcet and Q. Xu, \emph{Rosenthal
type inequalities for free chaos}. Preprint 2005.
\bibitem {JX} M. Junge and Q. Xu, \emph{Noncommutative
Burkholder/Rosenthal inequalities}. Ann. Probab. \textbf{31}
(2003), 948-995.
\bibitem {JX2} M. Junge and Q. Xu, \emph{On the best constants in
some non-commutative martingale inequalities}. Bull. London Math.
Soc. \textbf{37} (2005), 243-253.
\bibitem {LuP} F. Lust-Piquard and G. Pisier,
\emph{Non-commutative Khintchine and Paley inequalities}. Ark.
Mat. \textbf{29} (1991), 241-260.
\bibitem {PR} J. Parcet and N. Randrianantoanina, \emph{Gundy's
decomposition  for non-commutative martingales and applications}.
Preprint 2004.
\bibitem {PX1} G. Pisier and Q. Xu, \emph{Non-commutative
martingale inequalities}. Comm. Math. Phys. \textbf{189} (1997),
667-698.
\bibitem {PX2} G. Pisier and Q. Xu, \emph{Non-commutative
$L_p$-spaces}. Handbook of the Geometry of Banach Spaces II (Ed.
W.B. Johnson and J. Lindenstrauss) North-Holland (2003),
1459-1517.
\bibitem {R} N. Randrianantoanina, \emph{Non-commutative
martingale transforms}. J. Funct. Anal. \textbf{194} (2002),
181-212.
\bibitem {R2} N. Randrianantoanina, \emph{Square function
inequalities for non-commutative martingales}. Israel J. Math.
\textbf{140} (2004), 333-365.
\bibitem {R4} N. Randrianantoanina, \emph{Conditioned square
functions for non-commutative martingales}. Preprint 2005.
\bibitem {Ro} H.P. Rosenthal, \emph{On the subspaces of $L^p$
$(p>2)$ spanned by sequences of independent random variables}.
Israel J. Math. \textbf{8} (1970), 273-303.
\bibitem {W} F. Weisz, Martingale Hardy Spaces and its
Applications in Fourier Analysis. Lecture Notes in Math.
\textbf{1568}, Springer-Verlag, 1994.
\bibitem {X} Q. Xu,
\emph{Recent devepolment on non-commutative martingale
inequalities}. Functional Space Theory and its Applications.
Proceedings of International Conference \& 13th Academic Symposium
in China. Ed. Research Information Ltd UK. Wuhan 2003, 283-314.
\end{thebibliography}

\end{document}